\documentclass[11pt]{article} 
\usepackage[utf8]{inputenc} 
\usepackage[T1]{fontenc} 
\usepackage{comment} 
\usepackage{authblk} 
\usepackage{relsize} 
\usepackage{changepage} 
\usepackage{enumitem}
\usepackage{mathtools} 
\usepackage[english]{babel} 
\usepackage{amsmath} 
\usepackage{amsthm} 
\usepackage{amsfonts} 
\usepackage{url}
\usepackage{amssymb} 
\usepackage{graphicx} 
\usepackage{tabularx}
\usepackage[usenames,dvipsnames]{color} 
\theoremstyle{plain} 
\newtheorem{theorem}{Theorem}[section]
\newtheorem{corollary}[theorem]{Corollary} 
 
\newtheorem{lemma}[theorem]{Lemma}

\makeatletter 
\newcommand{\vast}{\bBigg@{4}} 
\newcommand{\Vast}{\bBigg@{5}} 
\makeatother 
\definecolor{bulgarianrose}{rgb}{0.28, 0.02, 0.03} 
\definecolor{gray}{rgb}{0.5, 0.5, 0.5} 
 
\theoremstyle{definition}

\theoremstyle{remark} 
\newtheorem*{remark}{Remark} 
\usepackage[left=2cm,right=2cm,top=2cm,bottom=2cm]{geometry} 
\usepackage{amssymb} 
\makeatletter 
\def\namedlabel#1#2{\begingroup
    #2%
    \def\@currentlabel{#2}%
    \phantomsection\label{#1}\endgroup
} 
\usepackage[normalem]{ulem} 
\usepackage{dsfont} 
\usepackage{accents} 
\usepackage{verbatim} 
\usepackage{ulem} 
\usepackage{pgf,tikz,pgfplots} 
\pgfplotsset{compat = 1.16} 
\usepackage[all]{xy} 
 
\usepackage{stackengine} 
\stackMath 
\newcommand\tsup[2][2]{%
 \def\useanchorwidth{T}%
  \ifnum#1>1%
    \stackon[-.5pt]{\tsup[\numexpr#1-1\relax]{#2}}{\scriptscriptstyle\sim}%
  \else%
    \stackon[.5pt]{#2}{\scriptscriptstyle\sim}%
  \fi%
}

\usepackage{enumerate} 
 
\usepackage{amsmath} 
\usepackage{amsfonts} 
\usepackage{amsthm} 
\usepackage{amssymb} 
\usepackage{array} 
\usepackage{booktabs} 
\usepackage{color} 
\usepackage{changepage} 
\usepackage[english]{babel} 
\usepackage{enumerate} 
\usepackage{epsfig} 
\usepackage{float} 
\usepackage{framed} 
\usepackage{gensymb} 
\usepackage{graphicx} 
\usepackage[colorlinks=true, allcolors=.]{hyperref}
\usepackage{indentfirst} 
\usepackage{longtable} 
\usepackage{mathtools} 
\usepackage{multirow} 
\usepackage[normalem]{ulem} 
\usepackage{setspace} 
\usepackage{subcaption} 
\usepackage{textcomp} 
\usepackage[utf8]{inputenc} 
\usepackage{verbatim} 
\usepackage{enumerate} 
\usepackage{xlop} 
\usepackage{listings} 
\usepackage[none]{hyphenat} 

\newcommand{\MYhref}[3][blue]{\href{#2}{\color{#1}{#3}}}%

\newcommand{\Mod}[1]{\ (\mathrm{mod}\ #1)}

\def\and{%
  \end{tabular}%
  \hskip 1em \@plus.17fil%
  \begin{tabular}[t]{c}}%
  
\title{\scshape Proof of the Complete Presence of a Modulo $4$ Bias for the Semiprimes}  
 
\author{Miroslav Marinov, Nikola Gyulev\footnote{m.marinov1617@gmail.com, gyulev2718@gmail.com} } 
\date{} 
\begin{document} 

\maketitle 
 
\begin{abstract}  
    In 2016, Dummit, Granville, and Kisilevsky showed that the proportion of semiprimes (products of two primes) not exceeding a given $x$, whose factors are congruent to $3$ modulo $4$, is more than a quarter when $x$ is sufficiently large. They have also conjectured that this holds from the very beginning, that is, for all $x \geq 9$. Here we give a proof of this conjecture. For $x\geq 1.1 \cdot 10^{13}$ we take an explicit approach based on their work. We rely on classical estimates for prime counting functions, as well as on very recent explicit improvements by Bennett, Martin, O'Bryant, and Rechnitzer, which have wide applications in essentially any setting involving estimations of sums over primes in arithmetic progressions. All $x < 1.1 \cdot 10^{13}$ are covered by a computed assisted verification.
 \end{abstract} 
 
\section{Introduction}   
There are many natural questions regarding the structure of prime numbers in arithmetic progression. Throughout we denote by $\pi(x)$ the number of primes not exceeding $x$ and by $\pi(x;q,a)$ the number of primes not exceeding $x$ which are congruent to $a$ modulo $q$. The prime number theorem for arithmetic progressions \cite{PNTAP} gives 
\begin{equation*} 
\pi(x;4,3) \sim \pi(x;4,1) \sim \frac{x}{2\log x}.   
\end{equation*} 
(Throughout we write $f(x) \sim g(x)$ if $g(x) \neq 0$ for sufficiently large $x$ and $\lim_{x\to \infty} \frac{f(x)}{g(x)} = 1$.) Hence it makes sense to conjecture that the two counts will be close throughout or at least alternate somewhat symmetrically. However, it seems that $\pi(x; 4, 3) > \pi(x; 4, 1)$ much more often. This is the so-called \textit{Chebyshev's Bias}, which he first observed in $1853$, describing it in a letter to Fuss \cite{de2022carrera}. This phenomenon has been thoroughly investigated in the twentieth century \cite{hardy1916contributions, kaczorowski1993contribution, knapowski1962comparative, LittlewoodSignChang}. In $1994$ Rubinstein and Sarnak \cite{rubinstein1994chebyshev} proved that, conditionally on the Generalized Riemann hypothesis and the Grand simplicity hypothesis, when measuring with logarithmic density Chebyshev is right approximately $99.59\%$ of the time. In general, it turns out that $\pi(x;n,a) > \pi(x;n,b)$ more often, for any quadratic non-residue $a$ and quadratic residue $b$ with $\gcd(a,n) = \gcd(b,n) = 1$.

Chebyshev's bias naturally leads to similar questions where the main object is slightly different. In particular, we can consider the products of two primes (also called \textit{semiprimes}) less than a given $x$. 
This has been explored by Dummit, Granville and Kisilevsky \cite{dummit2016big} in $2016$, who suggest strong dominance from remainder $3$. The behaviour for general modulus again turns out to depend on comparing quadratic versus non-quadratic residues. 
\begin{theorem}[Dummit-Granville-Kisilevsky, 2016] 
\label{DGKTheorem} 
Let $\chi$ be a real non-principal Dirichlet character with conductor $d$. For $ \eta \in \{ -1,1 \} $ we have 
\vspace{-0.3em} 
\begin{equation*} 
\frac{ \#\{ m \leq x : m = pq \ \text{with} \ \chi (p) = \chi (q) = \eta \} }{ \frac{ 1 }{ 4 } \#\{ m \leq x : m = pq, (m,d) = 1 \} } = 1 + \eta \frac{ \mathcal{L}_\chi + o(1) }{ \log\log x }, \hspace{0.3em} where \hspace{0.3em} \mathcal{L}_{\chi} := \sum_{p} \frac{\chi(p)}{p} \neq 0            
\end{equation*} 
and $p$, $q$ represent $($not necessarily distinct$)$ primes. 
\end{theorem} 
\noindent (Throughout, for functions $f$ and $g$ such that $g(x) \neq 0$ for sufficiently large $x$, we write $f(x) = o(g(x))$ if $\lim_{x\to \infty} f(x)/g(x) = 0$, and $f(x) = O(g(x))$ if $|f(x)| \leq Cg(x)$ for some constant $C$.)  
It follows that the sign of $\mathcal{L}_\chi$ determines the outcome of the bias. For modulo 4 with the character $\chi_{-4}(p) = (-1)^{\frac{p-1}{2}}$ we have $\mathcal{L}_{\chi_{\tiny{-4}}}\approx -0.334$, so for large $x$ the semiprimes $pq \leq x $ with       
\begin{itemize} 
\item $ p \equiv q \equiv \alpha \Mod 4 $ are approximately $\displaystyle \frac{1}{4} + \chi_{-4}(a)\frac{0.0835}{\log\log x} $ of all. 
\item $p \equiv 1 \Mod 4 , \ q \equiv 3 \Mod 4 $ are approximately $\frac{1}{2}$ of all. 
\end{itemize} 
An important note is that in this bias from some point on there is full dominance, differing substantially from Chebyshev's where a change of the lead occurs infinitely often. Hence we expect that Chebyshev's bias cannot be used to prove results for products of two primes.  
Note also that Theorem \ref{DGKTheorem} is not enough on its own to justify that the bias for a general modulus $d$ is present, as we also need that for any Dirichlet character $\chi$ we have $\mathcal{L}_{\chi}\neq 0,$ which is currently open and also mentioned by Dummit, Granville, and Kisilevsky. A further extension of Theorem \ref{DGKTheorem} has been done by Hough \cite{hough2017lower}, who showed that the right-hand side is at least $\displaystyle 1 + \frac{\log\log\log x + O(1)}{\log\log x}$ for at least one $d$ with $d\leq x$. The generalization to products of more than two primes has been explored by Meng \cite{meng2017distribution}. Finally, Dummit, Granville, and Kisilevsky conjecture that the bias modulo $4$ is present from the very beginning. Our main result is a proof of this conjecture.   
\begin{theorem} 
\label{thm:main} 
For all $x \geq 9$ we have 
\[ \frac{\#\{m\leq x: m=pq \mbox{ with } \ p\equiv q \equiv 3 \Mod 4\}}{\#\{m\leq x: m = pq \mbox{ with } p, q \mbox{ odd}\}} > \frac{1}{4}.  \] 
\end{theorem} 
\noindent We adapt Dummit, Granville, and Kisilevsky's main ideas to statements with explicit bounds, obtaining the result for $x \geq 1.1 \cdot 10^{13}$. The proof for $x \leq 1.1 \cdot 10^{13}$ shall be computer assisted, and the procedure is described separately at the end.

\section{Auxiliary results on counting primes}  
   
We now introduce explicit bounds concerning prime counting functions. Firstly, we state an explicit version of the prime number theorem. 
  
\begin{theorem} [Dusart  \cite{dusart1998autour}, Rosser-Schoenfeld \cite{rosser1962approximate}] 
\label{theorem:PNT}
We have
\begin{equation} 
\frac{x}{\log x} \leq \pi(x)     \leq       \frac{x}{\log x}+\frac{x}{( \log x )^2}+\frac{2.51x}{( \log x )^3}     \\\nonumber     
\end{equation} 
where the lower bound holds for $ x \geq 17 $, and the upper bound is for $ x \geq 355991 $.  
\end{theorem}     
\noindent A corresponding version of the prime number theorem for arithmetic progressions is as follows. 
\begin{theorem}[Bennett, Martin, O'Bryant, Rechnitzer, Corollary 1.6 in \cite{bennett2018explicit} for $q=4$] For $ x \geq 800 $ and $ a \in \{ 1,3 \} $ we have  
\label{APPNT} 
\begin{equation} 
\frac{x}{2\log x} < \pi(x;4,a) < \frac{x}{2\log x}+\frac{5x}{4( \log x )^2} .  \\\nonumber 
\end{equation} 
\end{theorem} 

We will also use a version with better error term, but with larger lower bound for $x$. As usual, denote Li$\displaystyle (x) := \int_2^x \frac{\mathrm{d}t}{\log t}$ for $x\geq 2$.

\begin{theorem}
    [Bennett, Martin, O'Bryant, Rechnitzer, \cite{bennett2018explicit}, Corollary $1.7$ for $q=4$] For $ x \geq 10^6$ and $ a \in \{ 1,3 \} $ we have  
\label{APPNT:stronger} 
\begin{equation} 
\left|\pi(x;4,a) - \frac{1}{2}{\normalfont \mbox{Li}}(x)\right| < 0.027\frac{x}{\log^2x}.  \\\nonumber 
\end{equation} 
\end{theorem} 

\noindent Note that there is also a version for $x\geq 5438260589$, but this would require us to restrict our main proof to $x\geq 5438260589^2 \approx 2.957 \cdot 10^{19}$, which is too large for our purposes. 

\smallskip

\noindent Appropriate bounds for the logarithmic integral are as follows:

\begin{lemma} 
\label{LI} 
     We have
     $$ \frac{x}{\log x} + \frac{x}{\log^2 x} < {\normalfont \mbox{Li}}(x) < \frac{x}{\log x} + \frac{6x}{5\log^2 x}$$
     where the lower bound is for all $x\geq 190$ and the upper bound is for all $x\geq 10^6$.
\end{lemma}

\begin{proof}
    The lower bound comes from \cite[Lemma 5.8]{bennett2018explicit}. For the upper bound, proceed as in \cite[Lemma 5.9]{bennett2018explicit}, computing $h(10^6) < \frac{6}{5}$, where $h(x) = ($Li$(x) - \frac{x}{\log x})/ \frac{x}{\log^2 x}$. 
\end{proof}


\begin{corollary}
    \label{APPNT:simplified-strong}
    For $x\geq 10^6$ and $a\in \{1,3\}$ we have
    \[ \frac{x}{2\log x} + 0.473\frac{x}{(\log x)^2} < \pi(x;4,a) < \frac{x}{2\log x} + 0.627\frac{x}{\log^2 x}.\]
\end{corollary}

\smallskip

\noindent Because of applications of partial summation with the sequence with $1/n $ for prime $n$, we make use of the following bound.   
\begin{theorem}[Mertens \cite{mertens1874beitrag}] For all $ x \geq 3 $ we have 
\label{MertennsT} 
\begin{equation*} 
\left| \sum_{p\leq x} \frac{1}{p} - \log\log x - M \right| \leq \frac{4}{\log (x + 1)} + \frac{2}{x \log x} 
\end{equation*} 
where $M \approx 0.26149$ is the Meissel-Mertens constant.           
\end{theorem} 

\begin{remark}
    We shall also make use of the simpler version $\left|\sum_{p\leq x}\frac{1}{p} - \log \log x\right| \leq 1$ for $x\geq 3$. It holds for $x\geq 227$ since $0\leq M \leq 1$ and $\frac{4}{\log (x+1)} + \frac{2}{x\log x} \leq 1-M$ (the left-hand side is decreasing). Note that $\log \log x$ is in $[0.09, 0.48)$ for $x\in [3,5)$, in $[0.48, 1.05)$ for $x\in [5,17)$ and in $[1.05, 1.7]$ for $x\in [17, 227]$; while $\sum_{p\leq x} \frac{1}{p}$ is $\frac{5}{6}$ for $x\in [3,5)$, in $[0.83, 1.35)$ for $x\in [5,17)$ and in $[1.35, 1.97]$ for $x\in [17, 227]$, so indeed the inequality holds for all $x \in [3,227]$, as well.
\end{remark}

\smallskip 
            
\noindent We also need a bound for $\pi(x;4, 3)$, which includes small $x$.
 
\begin{corollary} 
\label{c2.5}
    We have $$\displaystyle \frac{x}{2\log x} < \pi(x;4,3) \leq \frac{x}{2\log x} + \frac{5x}{4(\log x)^2} $$
    where the lower bound is for all $x\geq 19$ and the upper bound is for all $x>0$.
\end{corollary} 

\begin{proof} 
    By Corollary \ref{APPNT:simplified-strong}, it suffices to justify it for all $x \leq 10^6$. We illustrate it only for the lower bound, as the proof for the upper bound follows the same strategy. For all integers in $\left[19, 10^6\right]$, a Mathematica computation shows it is correct. Regarding non-integer $x$, in each interval $(a, a+1)$, $a\in \mathbb{Z}$, the function 
    \[ h(x) = \frac{x}{2\log x} - \pi(x; 4, 3), \mbox{ with } h'(x) = \frac{\log x - 1}{2(\log x)^2} \]
    is differentiable, with positive derivative for $x\geq 19$. Hence $h(x)$ is monotonic in each $(a, a+1)$, so knowing that $h(x) > 0$ when $x$ is an integer completes the proof.
\end{proof} 
 
\section{Main proof} 
We are now ready to prove the main theorem for $x\geq 1.1 \cdot 10^{13}$.
\begin{proof}
[Proof of Theorem \ref{thm:main}] We count the products of two primes $pq \leq x$ with $p\equiv q \equiv 3 \pmod 4$. Any such product can be represented with $p \leq q \leq \frac{x}{p}$ where $p \leq \sqrt{x}$. The count for products $pq \leq x$ with $p$, $q$ odd is analogous. Hence the desired statement
\begin{equation*} 
\frac{\#\{m\leq x: m=pq \mbox{ with } \ p\equiv q \equiv 3 \Mod 4\}}{\#\{m\leq x: m = pq \mbox{ with } p, q \mbox{ odd}\}} > \frac{1}{4}            
\end{equation*} 

\noindent is equivalent to

\begin{equation*}           
    W(x) := \sum_{\substack{p\leq \sqrt{x} \\ p\equiv 3 \, 
 (4)}}\sum_{\substack{p\leq q \leq x/p \\ q\equiv 3 \, (4)}} 1 \, - \, \frac{ 1 }{ 4 }  \sum_{2 < p\leq \sqrt{x}}\sum_{2<q\leq x/p}1 > 0. 
\end{equation*}        

\noindent As a consequence of Corollary \ref{c2.5}, (since $x\geq 19^2$, we have $\frac{x}{p} \geq \sqrt{x} \geq 19$), we have 
 		$$ \sum_{\substack{p\leq q \leq x/p \\ q\equiv 3 \, (4)}} 1 = \sum_{\substack{q \leq x/p \\ q\equiv 3 \, (4)}} 1 - \sum_{\substack{q < p \\ q\equiv 3 \, (4)}} 1 \geq \frac{x}{2p\log \frac{x}{p}}  - \sum_{\substack{q < p \\ q\equiv 3 \, (4)}}1.$$  

\noindent Next, we have $$\sum_{\substack{q < p \\ q\equiv 3 \, (4)}} 1 \leq \frac{p}{2\log p} + \frac{5p}{4(\log p)^2} \hspace{0.3em}\mbox{  and  } \sum_{\substack{q < p \\ q\equiv 3 \, (4)}} 1 \leq \frac{p}{2\log p} + 0.627\frac{p}{(\log p)^2}$$ from Corollary $\ref{c2.5}$ for $p<10^6$ and Corollary \ref{APPNT:simplified-strong} for $p>10^6$, respectively. Now note the simple observation that $\frac{1 - \chi_{-4}(p)}{2}$ equals $1$ if $p\equiv 3 \Mod 4$ and $0$ if $p\equiv 1 \Mod 4$. Thus, for $x\geq 10^{12}$

\begin{align*}
    \sum_{\substack{p\leq \sqrt{x} \\ p\equiv 3 \, 
 (4)}}\sum_{\substack{p\leq q \leq x/p \\ q\equiv 3 \, (4)}} 1 & \geq \sum_{2 < p \leq \sqrt{x}} \bigg(\frac{1 - \chi_{-4}(p)}{2} \frac{x}{2p\log \frac{x}{p}}\bigg) - \sum_{\substack{p \leq \sqrt{x} \\ p\equiv 3 \, (4)}}\bigg(\frac{p}{2\log p} + 0.627\frac{p}{(\log p)^2}\bigg) - 0.623\sum_{\substack{p < 10^6 \\ p\equiv 3 \, (4)}}\frac{p}{(\log p)^2} \\ 
 & \geq \sum_{2 < p \leq \sqrt{x}} \bigg(\frac{1 - \chi_{-4}(p)}{2} \frac{x}{2p\log \frac{x}{p}}\bigg) - \sum_{\substack{p \leq \sqrt{x} \\ p\equiv 3 \, (4)}}\bigg(\frac{p}{2\log p} + 0.627\frac{p}{(\log p)^2}\bigg) -  6.667\cdot10^7
\end{align*}

\noindent with the last expression being equal to   
$$ \frac{1}{4}\sum_{2<p\leq \sqrt{x}} \frac{x}{p\log \frac{x}{p}} - \frac{x}{4\log x}\sum_{p>2}\frac{\chi_{-4}(p)}{p} + \frac{x}{4\log x}\sum_{p>\sqrt{x}}\frac{\chi_{-4}(p)}{p} $$ $$ - \frac{x}{4\log x} \sum_{2 < p \leq \sqrt{x}} \frac{\chi_{-4}(p)\log p}{p\log \frac{x}{p}} - \sum_{\substack{p\leq \sqrt{x} \\ p\equiv 3 (4)}}\bigg(\frac{p}{2\log p} + 0.627\frac{p}{(\log p)^2}\bigg) - 6.667\cdot10^7. $$ 
Next, by Theorem \ref{theorem:PNT} for $x\geq 1.3 \cdot 10^{11} > 355991^2$
$$     
		    	 \sum_{2 < p\leq \sqrt{x}}\sum_{2<q\leq x/p}1 < \sum_{2<p\leq \sqrt{x}}\sum_{q\leq x/p}1 \leq \sum_{2<p\leq \sqrt{x}} \frac{x}{p\log \frac{x}{p}} + \sum_{2 < p\leq \sqrt{x}}\frac{x}{p(\log \frac{x}{p})^2} + 2.51\sum_{2 < p\leq \sqrt{x}}\frac{x}{p(\log \frac{x}{p})^3}. 
$$       
Hence the main target $W(x)$ is bounded below by   
$$
     - \frac{x}{4\log x}\sum_{p>2}\frac{\chi_{-4}(p)}{p} + \frac{x}{4\log x}\sum_{p>\sqrt{x}}\frac{\chi_{-4}(p)}{p} - \frac{x}{4\log x} \sum_{2<p \leq \sqrt{x}} \frac{\chi_{-4}(p)\log p}{p\log \frac{x}{p}} $$ 
    $$ - \sum_{ \substack{ p\leq \sqrt{x} \\ p\equiv 3  (4) } }\bigg(\frac{p}{2\log p} + 0.627\frac{p}{(\log p)^2}\bigg)  -  \frac{ x }{ 4 }\sum_{2 < p\leq \sqrt{x}}\frac{1}{p(\log \frac{x}{p})^2}- 0.6275x\sum_{ 2<p\leq \sqrt{x} } \frac{1}{p ( \log \frac{x}{p} )^3 }-6.667\cdot10^7.  
$$

\noindent Now we deal with each of the six sums separately.

\begin{itemize}[leftmargin=*]     
  \item Regarding $\sum_{p>2}\frac{\chi_{-4}(p)}{p}$, the work of Dummit-Granville-Kisilevsky \cite{dummit2016big} asserts that the decimal expansion of this sum is $(-0.334\ldots)$. 
  
  \item For $\sum_{p>\sqrt{x}}\frac{\chi_{-4}(p)}{p}$ we use partial summation with the function $\frac{1}{t}$ and the sequence $\chi_{-4}(n)$ when $n \neq 2$ is prime and $0$ otherwise. This gives 
  \begin{align*} 
    \sum_{p>\sqrt{x}}\frac{\chi_{-4}(p)}{p} = 
        \frac{\pi(\sqrt{x};4,3)-\pi(\sqrt{x};4,1)}{\sqrt{x}} + \int_{\sqrt{ x }}^{ \infty } \frac{ \pi(t;4,1)     -     \pi(t;4,3) }{ t^2 }\mathrm{d}t. 
    \end{align*} 
    By Corollary \ref{APPNT:simplified-strong} we have for $t\geq 10^6$  
\begin{equation} 
|\pi(t;4,3) - \pi(t;4,1)|  \leq  0.154\frac{t}{(\log t)^2}. \\\nonumber           
\end{equation} 
and so for $x\geq10^{12}$ the integral converges, yielding the bound
\begin{equation} 
\sum_{p>\sqrt{x}}\frac{\chi_{-4}(p)}{p} \geq - \frac{0.616}{(\log x)^2} - 0.154\int_{\sqrt{x}}^{\infty} \frac{\mathrm{d}t}{t(\log t)^2} 
=  - \frac{ 0.308 }{  \log  x } - \frac{0.616}{(\log x)^2} .  \nonumber 
\end{equation}

\item Using partial summation and Corollary \ref{APPNT:simplified-strong} for $x\geq 10^{12}$ 
		\begin{align*} 
		    \sum_{\substack{p\leq \sqrt{x} \\ p\equiv 3 \, (  4  )}} \frac{p}{\log p} & =\sum_{\substack{p< 10^6\\ p\equiv 3 \, (  4  )}} \frac{p}{\log p} + \pi(\sqrt{x};4,3) \frac{2\sqrt{x}}{\log x} - \pi(10^6; 4, 3)\frac{10^6}{6 \log 10} - \int_{10^6}^{\sqrt{x}}\pi(t;4,3)\frac{\log t - 1}{(\log t)^2}\mathrm{d}t \\ 
		    & < -1.428\cdot10^9 + \frac{2x}{(\log x)^2} + \frac{5.016x}{(\log x)^3} - \int_{10^6}^{\sqrt{x}} \bigg(  \frac{ ( \log t - 1 )  t }{ 2( \log t )^3 }+0.473\frac{(\log t - 1)t}{(\log t)^4}\bigg) \mathrm{d}t \\ 
		    & = -1.428\cdot10^9+\frac{x}{(\log x)^2} + \frac{5.016x}{(\log x)^3} + \frac{10^{12}}{144(\log 10)^2}\\& - 0.473\bigg[\frac{2}{3}\mbox{Li}(t^2)-\frac{t^2}{3\log t}-\frac{t^2}{6(\log t)^2}+\frac{t^2}{3(\log t)^3}\bigg]_{t=10^6}^{\sqrt{x}} 
            \end{align*} 
            as the antiderivative of $(\log t - 1)t/(2(\log t)^3)$ is $t^2/(2\log^2 t)$. Then, applying Lemma \ref{LI}, we see $\displaystyle \frac{2}{3}\mbox{Li}(t^2)-\frac{t^2}{3\log t}-\frac{t^2}{6(\log t)^2}>0$ for $t=\sqrt{x}$, thus obtaining the bound 
            \begin{align*}\sum_{\substack{p\leq \sqrt{x} \\ p\equiv 3 \, (  4  )}} \frac{p}{\log p} &  < \frac{x}{(\log x)^2}+\frac{3.755x}{(\log x)^3} - 2.4683 \cdot10^{7}.  
		\end{align*} 
{\small Similarly, we compute 
            \begin{align*}
		    \sum_{\substack{p\leq \sqrt{x} \\ p\equiv 3 \, (  4  )}} \frac{p}{(\log p)^2} &  =  \sum_{\substack{p< 10^6\\ p\equiv 3 \, (  4  )}} \frac{p}{(\log p)^2} + \pi(\sqrt{x}; 4, 3)\frac{4\sqrt{x}}{(\log x)^2} - \pi(10^6; 4,3) \frac{10^6}{36(\log 10)^2} - \int_{10^6}^{\sqrt{x}} \pi(t;4,3)\frac{\log t - 2}{(\log t)^3}\mathrm{d}t \\
            & < -9.9 \cdot 10^7 + \frac{4x}{(\log x)^3} + \frac{10.032x}{(\log x)^4} - \int_{10^6}^{\sqrt{x}} \left(\frac{(\log t - \frac{3}{2})t}{2(\log t)^4}+0.473\frac{(\log t - 2)t}{(\log t)^5}\right)\mathrm{d}t + \frac{1}{4}\int_{10^6}^{\sqrt{x}}\frac{t}{(\log t)^4}\mathrm{d}t
            \\ &=  -9.9\cdot10^7+\frac{4x}{(\log x)^3}+\frac{10.032x}{(\log x)^4} - \bigg[ \frac{t^2}{4(\log t)^3} \bigg]_{t=10^6}^{\sqrt{x}}  - 0.473\bigg[ \frac{t^2}{2(\log t)^4} \bigg]_{t=10^6}^{\sqrt{x}} + \frac{1}{4}\int_{10^6}^{\sqrt{x}}\frac{t}{(\log t)^4}\mathrm{d}t \\
            & < 2.3 \cdot 10^6 + \frac{2x}{(\log x)^3} + \frac{6.248x}{(\log x)^4} + \frac{1}{4}\int_{10^6}^{\sqrt{x}}\frac{t}{(\log t)^4} < 2.3 \cdot 10^6 + \frac{2x}{(\log x)^3} + \frac{10.248x}{(\log x)^4}
            \end{align*}  }where in the last step we used that $\frac{t}{(\log t)^4}$ is an increasing function for $t\geq e^4$ and hence bounded the integral from above by the maximum of the integrand times the length $\sqrt{x} - 10^6 < \sqrt{x}$ of the interval. Therefore 
\begin{equation*} 
\sum_{ \substack{ p\leq \sqrt{x} \\ p\equiv 3  (4) } }\bigg(\frac{p}{2\log p} + 0.627\frac{p}{(\log p)^2}\bigg)  \leq 
\frac{x}{2(\log x)^2} + \frac{3.1315x}{(\log x)^3} + \frac{6.43x}{(\log x)^4} - 1.08 \cdot 10^7.
\end{equation*} 

\item For a fixed $x$, we use partial summation with the function $1/\left(\log \frac{x}{t}\right)^2$ (on the variable $t$) and the sequence $\frac{1}{n}$ for $n$ prime. Together with Mertens' estimate in the form of Theorem \ref{MertennsT}, but with the right-hand side replaced by the slightly larger $\frac{81}{20\log x}$, we obtain
		\begin{align*}
		   \sum_{3 \leq p\leq \sqrt{x}} \frac{1}{p(\log \frac{x}{p})^2} & = \frac{\sum_{p\leq \sqrt{x}} \frac{1}{p} - \frac{1}{2}} {(\log \sqrt{x})^2} - \int_3^{\sqrt{x}}\left(\sum_{p\leq t} \frac{1}{p} - \frac{1}{2}\right)\frac{2}{t(\log \frac{x}{t})^3}\mathrm{d}t  \\   
		    & \leq \frac{4\left(\log \log \sqrt{x} - 0.23 + \frac{81}{20\log \sqrt{x}}\right)}{(\log x)^2} - \int_3^{\sqrt{x}}\left(\log \log t - \frac{81}{20\log t} - 0.24 \right)\frac{2}{t(\log \frac{x}{t})^3}\mathrm{d}t  \\ 
		    & = \frac{4\log \log x}{(\log x)^2} - \frac{4\log 2 + 0.92}{(\log x)^2} + \frac{32.4}{(\log x)^3} - \bigg[\frac{81\left(\log\log\frac{x}{t} - \log\log t\right)}{10\left(\log x\right)^3} + \frac{\log \log \frac{x}{t} - \log \log t}{(\log x)^2} \\ & + \frac{\log \log t}{\left(\log \frac{x}{t}\right)^2} - \frac{1}{\log x \log \frac{x}{t}} - \frac{81}{10(\log x)^2\left(\log \frac{x}{t}\right)} - \frac{81}{20(\log x)\left(\log \frac{x}{t}\right)^2} - \frac{6}{25\left(\log \frac{x}{t}\right)^2} \bigg]_{t=3}^{\sqrt{x}} \\ 
            & = \frac{2.04}{(\log x)^2} +\frac{64.8}{(\log x)^3} + A(x)
            \end{align*}
            where $A(x)$ is the exact antiderivative of $\left(\log \log t - \frac{81}{20\log t} - 0.24 \right)\frac{2}{t(\log \frac{x}{t})^3}$, evaluated at the lower limit $3$. In each occurrence of $\log \frac{x}{3}$ in a denominator, except for $\frac{\log \log 3}{\left(\log \frac{x}{3}\right)^2}$, we may replace it with $\log x$ since the corresponding term is negative. Hence
            \begin{align*}
                \sum_{3 \leq p\leq \sqrt{x}} \frac{1}{p(\log \frac{x}{p})^2} < \frac{0.8}{(\log x)^2} +\frac{52.65}{(\log x)^3} + \frac{81\left(\log\log\frac{x}{3} - \log\log 3\right)}{10\left(\log x\right)^3} + \frac{\log \log \frac{x}{3} - \log \log 3}{(\log x)^2} + \frac{\log \log 3}{\left(\log \frac{x}{3}\right)^2}. 
            \end{align*}

            For the last fraction we use $\frac{1}{\log \frac{x}{3}} \leq \frac{1.05}{\log x}$ for $x\geq 10^{12}$. Overall, we reach 
            \[ \sum_{3 \leq p\leq \sqrt{x}} \frac{1}{p(\log \frac{x}{p})^2} \leq \frac{\log \log \frac{x}{3}}{(\log x)^2} + \frac{0.81}{(\log x)^2} + \frac{8.1\log \log \frac{x}{3}}{(\log x)^3} + \frac{51.89}{(\log x)^3}.\]

\item An analogous approach, with the function $1/\left(\log \frac{x}{t}\right)^3$ (on the variable $t$),  yields 
\begin{align*} 
 \sum_{3 \leq p \leq \sqrt{x}} \frac{1}{p\left(\log\frac{x}{p}\right)^3} & =  \frac{\sum_{p\leq \sqrt{x}}\frac{1}{p} - \frac{1}{2}}{(\log\sqrt{x})^3}  - 3\int_3^{\sqrt{x}}\left(\log \log t - \frac{81}{20\log t} - 0.24 \right)\frac{1}{t(\log \frac{x}{t})^4}\mathrm{d}t \\ &
     \leq  \frac{ 8\left(\log\log \sqrt{x} - 0.23 + \frac{81}{20\log\sqrt{x}}\right)} { (\log x)^3 }  -  3 \int_{3}^{\sqrt{x}}     \frac{ \log \log t - \frac{81}{20\log t} - 0.24 }{ t ( \log \frac{x}{t})^4 }     \mathrm{d}t         \\   
    & \leq \frac{ 8\left(\log\log \sqrt{x} - 0.23 + \frac{81}{20\log\sqrt{x}}\right)} { (\log x)^3 }  -  \left(\frac{8\log \log \sqrt{x}}{(\log x)^3} - \frac{148}{25(\log x)^3} - \frac{81}{(\log x)^4}\right) + B(x) \\
    & \leq \frac{4.08}{(\log x)^3} + \frac{145.8}{(\log x)^4} + B(x)
\end{align*} 
where the value $B(x)$ of the integrand at the lower limit $3$ is
{\small 
\begin{align*}
    & \frac{\log \log \frac{x}{3} - \log \log 3}{(\log x)^3} + \frac{243}{20}\frac{\log \log \frac{x}{3} - \log \log 3}{(\log x)^4}  + \frac{\log \log 3 - \frac{87}{50}}{\left(\log \frac{x}{3}\right)^3} + \frac{\frac{5}{2}\log 3 - \frac{891}{40}}{\left(\log \frac{x}{3}\right)^3\log x} + \frac{\frac{243}{8}\log 3 - (\log 3)^2}{\left(\log \frac{x}{3}\right)^3(\log x)^2} - \frac{\frac{243}{20}(\log 3)^2}{\left(\log \frac{x}{3}\right)^3(\log x)^3} \\ & < \frac{\log \log \frac{x}{3}}{(\log x)^3} + \frac{243\log \log \frac{x}{3}}{20(\log x)^4} - \frac{1.74}{(\log x)^3} - \frac{20.68}{(\log x)^4} + \frac{36.33}{(\log x)^5} - \frac{14.67}{(\log x)^6}  
\end{align*} }
where the last inequality holds for $x\geq 10^{10}$. Overall, we reach
$$  \sum_{3 \leq p \leq \sqrt{x}} \frac{1}{p\left(\log\frac{x}{p}\right)^3} < \frac{\log \log \frac{x}{3}}{(\log x)^3} + \frac{12.15\log \log \frac{x}{3}}{(\log x)^4} + \frac{2.34}{(\log x)^3} + \frac{125.12}{(\log x)^4} + \frac{36.33}{(\log x)^5} - \frac{14.67}{(\log x)^6}. $$

\item For the sum $\displaystyle \sum_{3 \leq p \leq \sqrt{x}} \frac{\chi_{-4}(p)\log p}{p\log \frac{x}{p}}$, we split into two ranges for $p$ -- $[1,10^6)$ and $\left(10^6,\sqrt{x}\right]$. 

\vspace{-0.7em}
\begin{itemize}
    \item[$\blacktriangleright$] We derive a bound for $\sum_{p<10^6} \frac{\chi_{-4}(p)\log p}{p\log \frac{x}{p}}$ by generalizing an approach by Martin \cite{5031204}. It makes use of the following polynomial bounds, where $\ell \geq 0$ is an integer to be chosen later 
    $$\sum_{i=0}^{\ell} s^i < \frac{1}{1 - s} < s^{\ell} + \sum_{i=0}^{\ell} s^i, \text{ \hspace{0.3em}for } s \in \left(0, \frac{1}{2}\right).$$
    We apply it with $s=\frac{\log p}{\log x}\in(0,\frac{1}{2})$, for $x>10^{12}$ and $p<10^6$, which gives


$$ \sum_{i=0}^{\ell} \left(\frac{\log p}{\log x}\right)^{i} <\frac{1}{1-\frac{\log p}{\log x}}<  \left(\frac{\log p}{\log x}\right)^{\ell} + \sum_{i=0}^{\ell} \left(\frac{\log p}{\log x}\right)^{i}. $$
By writing $\displaystyle \frac{1}{\log \frac{x}{p}} = \frac{1}{\log x}\frac{1}{1-\frac{\log p}{\log x}}$, this leads to the upper bound

\begin{align*} 
&\frac{1}{\log x}\sum_{\substack{p<10^6\\p\equiv1(4)}}\left(\left(\frac{\log p}{\log x}\right)^{\ell} + \sum_{i=0}^{\ell} \left(\frac{\log p}{\log x}\right)^{i}\right)\frac{\log p}{p}- \frac{1}{\log x}\sum_{\substack{p<10^6\\p\equiv3(4)}}\left(\sum_{i=0}^{\ell} \left(\frac{\log p}{\log x}\right)^{i}\right)\frac{\log p}{p} \\ &= \sum_{i=1}^{\ell} \frac{1}{(\log x)^i}\sum_{p<10^6}\frac{\chi_{-4}(p)(\log p)^i}{p} +\frac{2}{(\log x)^{\ell+1}}\sum_{\substack{p<10^6\\p\equiv1(4)}}\frac{(\log p)^{\ell+1}}{p}-\frac{1}{(\log x)^{\ell+1}}\sum_{\substack{p<10^6\\p \equiv 3(4)}}\frac{(\log p)^{\ell+1}}{p}.
\end{align*} 
\noindent For $\ell = 6$ this is $\displaystyle -\frac{0.545}{\log x} - \frac{1.458}{(\log x)^2} - \frac{6.5}{(\log x)^3} - \frac{40.9}{(\log x)^4} - \frac{322}{(\log x)^5} + \frac{573393}{(\log x)^6}$.

\item[$\blacktriangleright$] The rest we split into residue classes and use partial summation (with $\frac{\log t}{t\log \frac{x}{t}}$) to get 
		\begin{align*}
		&\sum_{10^6< p \leq \sqrt{x}} \frac{\chi_{-4}(p)}{p\log \frac{x}{p}} = \sum_{b\in \{1,3\}}\chi_{-4}(b)\sum_{\substack{10^6 < p\leq\sqrt{x} \\ p\equiv b \ (4)}} \frac{\log p}{p\log \frac{x}{p}} \\
		& = \sum_{b \in \{1,3\}}\chi_{-4}(b)\left[\frac{1}{\sqrt{x}}\sum_{\substack{10^6< p\leq \sqrt{x} \\ p\equiv b \ (4)}} 1 - \int_{10^6}^{\sqrt{x}}\frac{\log x - \log t\log(\frac{x}{t})}{t^2(\log \frac{x}{t})^2} \sum_{\substack{10^6< p\leq t \\ p\equiv b \ (4)}}1 \mathrm{d}t\right].
		\end{align*}
		There are $78497$ odd primes less than $10^6$, from which $39322$ are $3$ modulo $4$ and $39175$ are $1$ modulo $4$. Hence 
        the above becomes
		$$\sum_{b\in \{1,3\}}\chi_{-4}(b)\left[\frac{\pi(\sqrt{x};4,b) - \pi(10^6;4,b)}{\sqrt{x}} + \int_{10^6}^{\sqrt{x}}\frac{\log t\log(\frac{x}{t})-\log x}{t^2(\log \frac{x}{t})^2} (\pi(t;4,b)-\pi(10^6;4,b)) \mathrm{d}t\right] $$
		$$ = \frac{\pi(\sqrt{x};4,1) - \pi(\sqrt{x};4,3) + 147}{\sqrt{x}} + \int_{10^6}^{\sqrt{x}}\frac{\log t\log(\frac{x}{t})-\log x}{t^2(\log \frac{x}{t})^2}(\pi(t;4,1)- \pi(t;4,3)+147)\mathrm{d}t.$$

       We now apply Corollary \ref{APPNT:simplified-strong}. Since $\log t \log \frac{x}{t} - \log x$ is positive for $3 \leq t \leq \sqrt{x}$ and $x\geq 3 \cdot 10^5$ (since $\log t \log \frac{x}{t}$ increases as $t$ increases), we obtain as an upper bound for $\sum_{3 \leq p \leq \sqrt{x}} \frac{\chi_{-4}(p)\log p}{p\log \frac{x}{p}}$ the expression 
		\begin{align*} 
	   \frac{147}{\sqrt{x}}+\frac{0.616}{(\log x)^2} +   
         \int_{10^6}^{\sqrt{x}}\frac{\log t\log(\frac{x}{t})-\log x}{t^2(\log \frac{x}{t})^2}\left(\frac{0.154t}{(\log t)^2}+147\right)\mathrm{d}t.
	\end{align*}  
\noindent Surprisingly, the contribution of the upper limit $\sqrt{x}$ from the integral is exactly $-\frac{147}{\sqrt{x}}$. Hence we reach as an upper bound  
{\footnotesize 
\begin{align*}
    &\frac{0.616}{(\log x)^2} + \frac{147\log 10^6}{10^6\log \frac{x}{10^6}} + 0.154 \left(\frac{\log \log \frac{x}{10^6} - \log \log 10^6}{\log x} - 2\frac{\log \log \frac{x}{10^6} - \log\log 10^6}{(\log x)^2} - \frac{1}{\log x\log 10^6} + \frac{2}{\log \frac{x}{10^6}\log x }\right) \\
    & \leq  \frac{0.154\log \log \frac{x}{10^6}}{\log x} - \frac{0.308\log \log \frac{x}{10^6}}{(\log x)^2} - \frac{0.415}{\log x} + \frac{0.0021}{\log \frac{x}{10^6}} +\frac{1.425}{(\log x)^2} + \frac{0.308}{\log \frac{x}{10^6}{\log x}}  
\end{align*}}
\end{itemize}
\noindent Overall, we have reached
{\small
\[ \sum_{3 \leq p \leq \sqrt{x}} \frac{\chi_{-4}(p)\log p}{p\log \frac{x}{p}} < \frac{0.154\log \log \frac{x}{10^6}}{\log x} - \frac{0.308\log \log \frac{x}{10^6}}{(\log x)^2} - \frac{0.955}{\log x} + \frac{0.583}{(\log x)^2}  - \frac{6.5}{(\log x)^3} - \frac{40.9}{(\log x)^4} - \frac{322}{(\log x)^5} + \frac{573393}{(\log x)^6}. \]}
\end{itemize} 
\noindent \textbf{Collecting all bounds}. In conclusion, the main target $W(x)$ is bounded below by

$$ \frac{0.0835x}{\log x}-\frac{0.54075x}{(\log x)^2}-\frac{17.8721x}{(\log x)^3}-\frac{83.3178x}{(\log x)^4}-\frac{12.5721x}{(\log x)^5}+\frac{89.7054x}{(\log x)^6}-\frac{143348.25x}{(\log x)^7}-\frac{0.0385x\log\log\frac{x}{10^6}}{(\log x)^2} $$ $$+\frac{0.077x\log\log\frac{x}{10^6}}{(\log x)^3} -\frac{0.25x\log\log\frac{x}{3}}{(\log x)^2}-\frac{2.6525x\log\log\frac{x}{3}}{(\log x)^3}-\frac{7.62412x\log\log\frac{x}{3}}{(\log x)^4}-5.587\cdot10^7 $$

\noindent for $x\geq 10^{12}$. This is positive for $x\geq 1.1 \cdot 10^{13}$ by a direct derivative computation. \end{proof} 

\section{Computer-assisted proof for the remaining $x$}

Computations with large primes generally require a lot of time and memory. For the sake of reducing the need of these, our verification will not give precise values of the ratio of products with two primes $3$ mod $4$ to all products. Most of the time, it will provide a lower bound for the numerator and an upper bound for the denominator. Hence, it gives a lower bound for the ratio, which turns out to always exceed $\frac{1}{4}$ but is always reasonably close to $\frac{1}{4}$. The code and output are available at a public \MYhref{https://github.com/Assassino9931/Prime4-code/blob/main/Prime4.c}{GitHub repository} of one of the authors. The most important value is {\tt max$\_$prime}, currently set to $10^9$ -- the full computation takes approximately 1 hour on a computer with an Intel 2.1 GHz processor and 8 GB RAM. The value of {\tt max$\_$prime} can be adjusted by the user to anything less than $10^{18}$ (due to int type limits in C), as long as it exceeds the square root of the desired goal -- in our case, the goal is $1.1 \cdot 10^{13}$.  
 
\smallskip 
 
\noindent \textbf{Step 1}. We begin with a useful precomputation. Let $P$ denote the product of all primes from $2$ to $23$. First, calculate all remainders modulo $P$, which are coprime to $P$. Then, on each prime up to {\tt max$\_$prime} we assign a boolean value $1$, using a Sieve of Eratosthenes and the previous information regarding numbers coprime with $P$. Afterwards, again up to {\tt max$\_$prime}, we compute and store the number of primes $3$ mod $4$ and the overall number of primes up to $2^{13} \cdot k$ for each integer $k$. Finally, we make a separate list of all primes up to $10^7$ (which is enough for the purposes of products bounded by $1.1 \cdot 10^{13}$) in order to access them faster in the next steps. In total, we have three lists, each one being a subset of the previous, which will be used for navigational purposes. 

\smallskip

\noindent \textbf{Step 2}. We implement a simple idea that significantly speeds up the process. Suppose we have bounded the numerator $n_1$ and denominator $d_1$ of the ratio for $x = x_1$ (initially $x=9$) and that the ratio indeed exceeds $\displaystyle \frac{1}{4}$. The main observation is that we can now move the computation to $x = x_1 + 2(n_1 - d_1)$, since for all $x$ in between the ratio still exceeds $\displaystyle \frac{1}{4}$, even if all odd integers in between are primes $1$ mod $4$. The factor $2$ comes for free since we skip even integers, as they are not prime. Note that in the process, we will not update and use previous values of $n_1$ and $d_1$, but only $x$. 

\smallskip 

\noindent \textbf{Step 3}. Finally, we do the actual computation. For each $x$ we compute the denominator as follows: we loop through all primes $p$ not exceeding $\sqrt{x}$ and then, depending on whether $\displaystyle \frac{x}{p} \leq $ {\tt max$\_$prime} or not, we either give an exact count for the number of possibilities for the second prime, or we bound it using Theorem \ref{theorem:PNT} and Corollary \ref{APPNT:simplified-strong}. When giving the exact count, we use the navigation between lists mentioned at the end of Step 1, as this speeds up the procedure. To avoid overcounting, we have to subtract the index of the current prime $p$ in the loop. The numerator is separately computed with primes $3$ mod $4$ in an analogous manner.

\smallskip

For each needed $x$, we have obtained a lower bound for the ratio, showing that it is always greater than $\displaystyle \frac{1}{4}$. The code is designed to stop if some $x$ does not satisfy the desired hypothesis, i.e., the obtained lower bound is not greater than $\displaystyle \frac{1}{4}$. In an actual run of the code, one sees there is no such stop, so the hypothesis is indeed correct. 

\bigskip 

\noindent \textbf{Acknowledgments.} We would like to thank James Maynard for discussions with one of the authors on the topic during his master's thesis work at the University of Oxford. We are indebted to Tsvetomir Kaydzhiev and Boris Mihov for essential technical assistance with the computer verification of the hypothesis for those $x$ where our methods were out of reach. We are grateful to Angel Raychev for proof-checking particular technical details, as well as Iliyas Noman, Radostin Cholakov, Vladimir Mitankin, and representatives of the American Mathematical Society and Mu Alpha Theta at Regeneron ISEF 2023 for general comments on the paper. Thanks also to a reviewer for a lot of technical corrections and advice on the exposition of the paper. This research was supported by the Summer Research School at the High School Student Institute of Mathematics and Informatics (HSSIMI) under an ``Education with Science'' grant in Bulgaria and the provided atmosphere for work is greatly appreciated.          

\bibliographystyle{plain} 
\bibliography{References}

\begin{thebibliography}{10}

\bibitem{bennett2018explicit}
M.~A.~Bennett, K.~O’Bryant, G.~Martin and A.~Rechnitzer.
\newblock Explicit bounds for primes in arithmetic progressions.
\newblock {\em Illinois Journal of Mathematics}, 62(1-4):427--532, 2018.

\bibitem{de2022carrera}
P.~L. Chebyshev.
\newblock Lettre de m. le professeur tchébychev à m. fuss sur un nouveaux
  théorème relatif aux nombres premiers contenus dans les formes 4n + 1 et 4n
  + 3.
\newblock {\em Bull. Phys. Acad. Sci. St. Petersburg}, 11, 1853.

\bibitem{PNTAP}
Ch. J. De~la Vallee~Poussin.
\newblock Recherches analytiques sur la théorie des nombres premiers.
\newblock {\em Annales de la Société scientifique de Bruxelles, Imprimeur de
  l'Académie Royale de Belgique}, 1896.
  
\bibitem{dummit2016big}
D.~Dummit, A.~Granville and H.~Kisilevsky.
\newblock Big biases amongst products of two primes.
\newblock {\em Mathematika}, 62(2):502--507, 2016.

\bibitem{dusart1998autour}
P.~Dusart.
\newblock Autour de la fonction qui compte le nombre de nombres premiers.
\newblock 1998.

\bibitem{hardy1916contributions}
G.~H. Hardy and J.~E. Littlewood.
\newblock Contributions to the theory of the riemann zeta-function and the
  theory of the distribution of primes.
\newblock {\em Acta Mathematica}, 41:119--196, 1916.

\bibitem{hough2017lower}
P.~Hough.
\newblock A lower bound for biases amongst products of two primes.
\newblock {\em Research in Number Theory}, 3:1--11, 2017.

\bibitem{kaczorowski1993contribution}
J.~Kaczorowski.
\newblock A contribution to the shanks-r{\'e}nyi race problem.
\newblock {\em The Quarterly Journal of Mathematics}, 44(4):451--458, 1993.

\bibitem{knapowski1962comparative}
S.~Knapowski and P.~Tur{\'a}n.
\newblock Comparative prime-number theory. i.
\newblock {\em Acta Mathematica Hungarica}, 13(3):299--314, 1962.

\bibitem{LittlewoodSignChang}
J.~Littlewood.
\newblock Sur la distribution des nombres premiers.
\newblock {\em Comptes Rendus}, 1914.

\bibitem{5031204}
G.~Martin.
\newblock $\sum_{p \leq 10^6}\frac{(-1)^{\frac{p-1}{2}}\log p}{p\log \frac{x}{p}}$ is negative for large $x$.
\newblock Mathematics Stack Exchange, \url{https://math.stackexchange.com/q/5031204} (version: 2025-02-04).

\bibitem{meng2017distribution}
X.~Meng.
\newblock The distribution of k-free numbers and integers with fixed number of
  prime factors.
\newblock 2017.

\bibitem{mertens1874beitrag}
F.~Mertens.
\newblock Ein beitrag zur analytischen zahlentheorie.
\newblock 1874.

\bibitem{rosser1962approximate}
J.~B. Rosser and L.~Schoenfeld.
\newblock Approximate formulas for some functions of prime numbers.
\newblock {\em Illinois Journal of Mathematics}, 6(1):64--94, 1962.

\bibitem{rubinstein1994chebyshev}
M.~Rubinstein and P.~Sarnak.
\newblock Chebyshev's bias.
\newblock {\em Experimental Mathematics}, 3(3):173--197, 1994.

\end{thebibliography}

\end{document}